\documentclass{article}
\usepackage[utf8]{inputenc}

\usepackage{amsmath,dsfont}
\allowdisplaybreaks[4]

\usepackage[all]{xy}
\usepackage{amssymb}
\usepackage{amsthm}
\usepackage{hyperref}
\hypersetup{colorlinks=true,linkcolor=blue,citecolor=red}
\usepackage{amsmath}
\usepackage{amscd}
\usepackage{verbatim}
\usepackage{eurosym}
\usepackage{float}
\usepackage{color}
\usepackage{array}
\usepackage{dcolumn}
\usepackage[mathscr]{eucal}
\usepackage[all]{xy}
\usepackage{hyperref}
\usepackage{tikz}
\usetikzlibrary{decorations.pathmorphing}
\usetikzlibrary{decorations.pathreplacing}
\usetikzlibrary{decorations.markings}
\usetikzlibrary{decorations.footprints}
\usetikzlibrary{decorations.shapes}
\usetikzlibrary{decorations.fractals}
\usetikzlibrary{arrows,automata,chains,shadows}
\usepackage{mathrsfs}
\usepackage{amsmath}
\usepackage{amssymb}
\usepackage{amsfonts,ifpdf}
\usepackage{graphicx}
\usepackage{times}
\usepackage{float}
\usepackage{epstopdf}
\usepackage{cite}
\usepackage{youngtab}
\usepackage{ytableau}
\ytableausetup{mathmode, boxsize=0.9em}

\newtheorem{thm}{Theorem}[section]

\newtheorem{pf}{proof}[section]
\theoremstyle{remark}

\makeatletter \@addtoreset{equation}{section} \makeatother
\makeindex 
\def\qed{\hfill \rule{4pt}{7pt}}

\begin{document}

\begin{center}

 {\Large \bf  A Cyclic Analogue of Stanley's Shuffle Theorem}

\end{center}

\begin{center}
{Kathy Q. Ji}$^{1}$ and {Dax T.X. Zhang}$^{2}$  \vskip 2mm

$^{1,2}$ Center for Applied Mathematics,  Tianjin University, Tianjin 300072, P.R. China\\[6pt]
   \vskip 2mm

    $^1$kathyji@tju.edu.cn and $^2$zhangtianxing6@tju.edu.cn
\end{center}

\vskip 6mm \noindent {\bf Abstract.}   We introduce the cyclic major index of a cycle permutation and give a bivariate analogue of  enumerative formula for the cyclic shuffles with a given cyclic descent number due to Adin, Gessel, Reiner and Roichman, which can be viewed as a cyclic analogue of Stanley's Shuffle Theorem. This gives an answer to a   question of Adin, Gessel, Reiner and Roichman, which has been posed by   Domagalski, Liang, Minnich, Sagan, Schmidt and
Sietsema again.



 \vskip 6mm

\section{Introduction}

The main theme of this note is to establish a cyclic analogue of Stanley's Shuffle Theorem. Recall that Stanley's Shuffle Theorem establishes an explicit expression for the generating function of the number of shufflings of two disjoint permutation $\sigma$ and $\pi$  with a given cyclic descent number   and a  given  major index.
 Here we adopt  some common notation and terminology on permutations as used in \cite[Chapter 1]{Stanley-2012}. We say that $\pi=\pi_1\pi_2\cdots \pi_n$ is a permutation of length $n$   if it is a sequence of $n$ distinct letters (not necessarily from $1$ to $n$). For example, $\pi=9\,2\,8\,10\,12\,3\,7$ is a permutation of length $7$. Let  $\mathcal{S}_n$ denote the set of all permutations of length $n$.

 Let $\pi \in \mathcal{S}_n$, we say that  $1\leq i \leq n-1$ is a  descent of $\pi$  if $\pi_i>\pi_{i+1}$.  The set of  descents of $\pi$  is called the descent set of $\pi$, denoted ${\rm Des}(\pi)$, viz.,
 \[{\rm Des}(\pi):=\{1\leq i \leq n-1:\pi_i>\pi_{i-1}\}.\]
 The number of its descents is called the descent number, denoted ${\rm des}(\pi)$, namely,
  \[{\rm des}(\pi):=\# {\rm Des}(\pi),\]
where the hash symbol $\# \mathcal{S}$ stands for the cardinality  of a set $\mathcal{S}$.
 The   major index of $\pi$, denoted ${\rm maj}(\pi)$,  is defined to be the sum of its  descents. To wit,
\[{\rm maj}(\pi):=\sum_{k \in {\rm Des}(\pi)}k.\]

Let $\sigma \in \mathcal{S}_n$ and $\pi \in \mathcal{S}_m$ be  disjoint permutations, that is, permutations with no letters in common. We say that $\alpha \in \mathcal{S}_{n+m}$ is a shuffle of $\sigma$ and $\pi$ if both $\sigma$ and $\pi$  are
subsequences of $\alpha$. The set of shuffles of  $\sigma$ and $\pi$  is denoted $\mathcal{S}(\sigma, \pi)$. For example,
\[\mathcal{S}(6\,3,1\,4)=\{6\,3\,1\,4, 6\,1\,3\,4, 6\,1\,4\,3, 1\,4\,6\,3, 1\,6\,3\,4, 1\,6\,4\,3\}.\]
Clearly, the number of permutations in $\mathcal{S}(\sigma, \pi) $ is ${m+n \choose n}$ for two disjoint permutations $\sigma \in \mathcal{S}_n$ and $\pi \in \mathcal{S}_m$.

Stanley's Shuffle Theorem states that
\begin{thm}\label{stanley}
Let $\sigma \in \mathcal{S}_m$ and $\pi \in \mathcal{S}_n$ be disjoint permutations, where ${\rm des}(\sigma)=r$  and  ${\rm des}(\pi)=s$.  Then
\begin{equation}
   \sum_{\alpha\in \mathcal{S}(\sigma,\pi) \atop {\rm des}(\alpha)=k}q^{{\rm maj}(\alpha)}= {m-r+s \brack k-r} {n-s+r  \brack  k-s}  q^{{\rm maj}(\sigma)+{\rm maj}(\pi)+(k-s)(k-r)}.
   \end{equation}

\end{thm}
Here
\[{n \brack m}=\frac{(1-q^n)(1-q^{n-1})\cdots (1-q^{n-m+1})}{(1-q^m)(1-q^{m-1})\cdots (1-q)}\]
is the Gaussian polynomial (also  called the $q$-binomial coefficient), see Andrews \cite[Chapter 1]{Andrews-1976}.

  Stanley \cite{Stanley-1972}  obtained the above expression in light of the $q$-Pfaff-Saalsch\"utz identity in his setting of $P$-partitions.  The bijective proofs of  Stanley's shuffle theorem have been given by Goulden
  \cite{Goulden-1985},   Stadler \cite{Stadler-1999}, Ji and Zhang \cite{Ji-Zhang-2022},  respectively.

\begin{figure}[h!]
\begin{center}
\begin{tikzpicture}[scale=0.5]
\draw[line width=1pt] (0, 0) circle (2);  
\draw[line width=1pt] (0, 0) circle (1);
\node at(90:2)[below=2pt]{6};
\node at(0:2)[left]{3};
\node at(-90:2)[above]{{\bf 1}};
\node at(180:2)[right]{{\bf 4}};
\node at(-90:2)[below=17pt]{[6\,3\,{\bf 1}\,{\bf 4}]};
\foreach \x in {90,0,-90,180}
\draw [line width=1pt](\x:1.85)--(\x:2.15);
\draw[-latex, line width=1pt] (-1,0) arc(-179: -180: 1);
\end{tikzpicture}
\quad \quad \quad \quad \quad
\begin{tikzpicture}[scale=0.5]
\draw[line width=1pt] (0, 0) circle (2);  
\draw[line width=1pt] (0, 0) circle (1);
\node at(90:2)[below=2pt]{6};
\node at(0:2)[left]{3};
\node at(-90:2)[above]{{\bf 4}};
\node at(180:2)[right]{{\bf 1}};
\node at(-90:2)[below=17pt]{[6\,3\,{\bf 4}\,{\bf 1}]};
\foreach \x in {90,0,-90,180}
\draw [line width=1pt](\x:1.85)--(\x:2.15);
\draw[-latex, line width=1pt] (-1,0) arc(-179: -180: 1);
\end{tikzpicture}
\quad \quad \quad \quad \quad
\begin{tikzpicture}[scale=0.5]
\draw[line width=1pt] (0, 0) circle (2);  
\draw[line width=1pt] (0, 0) circle (1);
\node at(90:2)[below=2pt]{6};
\node at(0:2)[left]{{\bf 1}};
\node at(-90:2)[above]{{\bf 4}};
\node at(180:2)[right]{3};
\node at(-90:2)[below=17pt]{[6\,{\bf 1}\,{\bf 4}\,3]};
\foreach \x in {90,0,-90,180}
\draw [line width=1pt](\x:1.85)--(\x:2.15);
\draw[-latex, line width=1pt] (-1,0) arc(-179: -180: 1);
\end{tikzpicture}
\\[15pt]
\begin{tikzpicture}[scale=0.5]
\draw[line width=1pt] (0, 0) circle (2);  
\draw[line width=1pt] (0, 0) circle (1);
\node at(90:2)[below=2pt]{6};
\node at(0:2)[left]{{\bf 4}};
\node at(-90:2)[above]{{\bf 1}};
\node at(180:2)[right]{3};
\node at(-90:2)[below=17pt]{[6\,{\bf 4}\,{\bf 1}\,3]};
\foreach \x in {90,0,-90,180}
\draw [line width=1pt](\x:1.85)--(\x:2.15);
\draw[-latex, line width=1pt] (-1,0) arc(-179: -180: 1);
\end{tikzpicture}
\quad \quad \quad \quad \quad
\begin{tikzpicture}[scale=0.5]
\draw[line width=1pt] (0, 0) circle (2);  
\draw[line width=1pt] (0, 0) circle (1);
\node at(90:2)[below=2pt]{6};
\node at(0:2)[left]{{\bf 1}};
\node at(-90:2)[above]{3};
\node at(180:2)[right]{{\bf 4}};
\node at(-90:2)[below=17pt]{[6\,{\bf 1}\,3\,{\bf 4}]};
\foreach \x in {90,0,-90,180}
\draw [line width=1pt](\x:1.85)--(\x:2.15);
\draw[-latex, line width=1pt] (-1,0) arc(-179: -180: 1);
\end{tikzpicture}
\quad \quad \quad \quad \quad
\begin{tikzpicture}[scale=0.5]
\draw[line width=1pt] (0, 0) circle (2);  
\draw[line width=1pt] (0, 0) circle (1);
\node at(90:2)[below=2pt]{6};
\node at(0:2)[left]{{\bf 4}};
\node at(-90:2)[above]{3};
\node at(180:2)[right]{{\bf 1}};
\node at(-90:2)[below=17pt]{[6\,{\bf 4}\,3\,{\bf 1}]};
\foreach \x in {90,0,-90,180}
\draw [line width=1pt](\x:1.85)--(\x:2.15);
\draw[-latex, line width=1pt] (-1,0) arc(-179: -180: 1);
\end{tikzpicture}
\end{center}
\caption{The circular representations  of cyclic shuffles of $[6\,3]$ and $[4\,1]$.}\label{fig}
\end{figure}
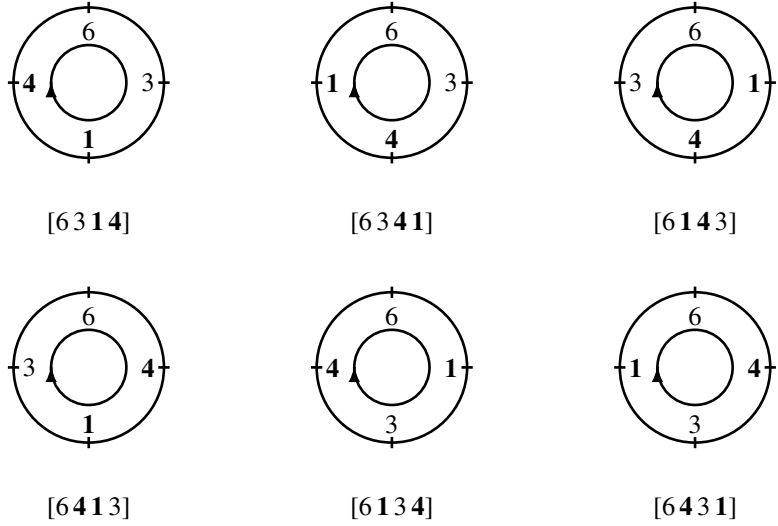

Recently, Adin, Gessel, Reiner and Roichman \cite{Adin-Gessel-Reiner-Roichman-2021}  introduced a cyclic version of quasisymmetric functions with a corresponding cyclic shuffle operation.  A cyclic permutation $[\pi]$ of length $n$ can be viewed as an equivalence class of linear permutations $\pi=\pi_1\pi_2\cdots \pi_n$ of length $n$ under the cyclic equivalence relation $\pi_1\pi_2\cdots \pi_n \sim \pi_i\cdots \pi_n \pi_1\cdots \pi_{i-1}$ for all $2\leq i\leq n$. For example,
\begin{equation}\label{examp-r}
[4\,2\,3\,1]=\{4\,2\,3\,1,2\,3\,1\,4,3\,1\,4\,2,1\,4\,2\,3\}
\end{equation}
is a cyclic permutation of length $4$, where
\[[4\,2\,3\,1]=[2\,3\,1\,4]=[3\,1\,4\,2]=[1\,4\,2\,3].\]
Let $\pi_l$ be the largest element in $[\pi]$,  the linear permutation $\hat{\pi}=\pi_l\pi_{l+1}\cdots \pi_n\pi_1\cdots \pi_{l-1}$ corresponding to the cyclic permutation $[\pi]$ is called the representative of the cyclic permutation $[\pi]$.
For the example above,   $4\,2\,3\,1$ is the representative of the cyclic permutation $[4\,2\,3\,1]$. Here and in the sequel, we use the representative to represent  each cyclic permutation $[\pi]$. For example, we use $[4\,2\,3\,1]$ to represent the equivalence class in \eqref{examp-r}.  In this way, all cyclic permutations of $\{1,2,3,4\}$ are listed as follows:
\[[4\,1\,2\,3],[4\,3\,1\,2], [4\,1\,3\,2], [4\,2\,1\,3],[4\,2\,3\,1],[4\,3\,2\,1].\]
Let  $c\mathcal{S}_n$ denote the set of all cyclic permutations of length $n$ and let $[\sigma] \in c\mathcal{S}_n$ and $[\pi] \in c\mathcal{S}_m$ be disjoint cyclic permutations, that is, cyclic permutations with no letters in common. We say that $[\alpha] \in c\mathcal{S}_{n+m}$ is a cyclic shuffle of two cyclic permutations $[\sigma]$ and $[\pi]$ if both $[\sigma]$ and $[\pi]$  are circular subsequences  of $[\alpha]$. The set of cyclic shuffles of $[\sigma]$ and $[\pi]$ is denoted $c\mathcal{S}([\sigma],[\pi])$. For example,
\begin{equation}\label{example-c}
c\mathcal{S}([6\,3],[4\,1])=\{[6\,3\,{\bf 1}\,{\bf 4}],[6\,3\,{\bf 4}\,{\bf 1}],[6\,{\bf 1}\,{\bf 4}\,3],[6\,{\bf 4}\,{\bf 1}\,3],[6\,{\bf 1}\,3\,{\bf 4}],[6\,{\bf 4}\,3\,{\bf 1}]\}.
\end{equation}
The elements of $[\pi]$ in $[\alpha]$ are in boldface to distinguish them from the elements of $[\sigma]$. Figure \ref{fig} lays out the circular representations of cyclic shuffles of $[6\,3]$ and $[4\,1]$.

Evidently,
\begin{equation}\label{cycli-enum}
   \# c\mathcal{S}([\sigma],[\pi])=(m+n-1){m+n-2\choose m-1},
\end{equation}
for two disjoint cyclic permutations $[\sigma] \in c\mathcal{S}_n$ and $[\pi] \in c\mathcal{S}_m$,  see \cite[Eq.~(7)]{Domagalski-Liang-Minnich-Sagan-Schmidt-Sietsema-2021}.

In order to study Solomon's descent algebra, Cellini  \cite{Cellini-1995, Cellini-1998} introduced  the cyclic descent set. Let $\pi=\pi_1\pi_2\ldots\pi_n$ be a linear permutation, the cyclic descent set of $\pi$ is defined to be
 \begin{equation*}
     {\rm cDes} (\pi)=\{1\leq i \leq n\colon \pi_i>\pi_{i+1}\}
     \end{equation*}
     with the convention $\pi_{n+1}=\pi_1$.  The number of its cyclic descents is called the cyclic descent number, denoted   $  {\rm cdes}(\pi)$, viz.,
 \[{\rm cdes}(\pi):=\# {\rm cDes}(\pi).\]
Let $[\pi]$ be a cyclic permutation of length $n$,  note that all linear permutations corresponding to $[\pi]$ have the same number of cyclic descents, so we may define the cyclic descent number of $[\pi]$ as
\begin{equation}\label{defi-cdesnum}
    {\rm cdes}\left([\pi]\right)={\rm cdes}\left(\pi\right),
\end{equation}
where $\pi$ is any linear permutation corresponding to $[\pi]$.

Based on their setting of cyclic quasi-symmetric functions,  Adin, Gessel, Reiner and Roichman \cite{Adin-Gessel-Reiner-Roichman-2021} established the following enumerative formula for the cyclic shuffles with a given cyclic descent number.

\begin{thm}[Adin-Gessel-Reiner-Roichman]\label{AGRR} Let $[\sigma]\in c\mathcal{S}_m$ and $[\pi] \in c\mathcal{S}_n$ be disjoint cyclic permutations, where ${\rm cdes}([\sigma])=r$ and ${\rm cdes}([\pi])=s$. Let $c\mathcal{S}([\sigma],[\pi],k)$ denote the set of cyclic shuffles of $[\sigma]$ and $[\pi]$ with $k$ cyclic descent number. Then
\begin{align}\label{AGRR-e}
  \# c\mathcal{S}([\sigma],[\pi],k)  = \frac{k(m-r)(n-s)+(m+n-k)rs}{(m-r+s)(n-s+r)}\dbinom{m-r+s}{k-r}\dbinom{n-s+r}{k-s}.
\end{align}
\end{thm}
Summing  \eqref{AGRR-e} over all $k$ gives \eqref{cycli-enum} upon  using the Chu-Vandermond identity \cite[Eq.~(3.3.10) with $q=1$]{Andrews-1976}. At the end of their paper,
Adin, Gessel, Reiner and Roichman \cite{Adin-Gessel-Reiner-Roichman-2021} asked a question  about looking for a notion of cyclic major index, which provides a bivariate analogue of Theorem \ref{AGRR}. This question has been posed  by Domagalski, Liang, Minnich, Sagan, Schmidt and Sietsema  in \cite[Question 4.1]{Domagalski-Liang-Minnich-Sagan-Schmidt-Sietsema-2021} again.

In this paper, we introduce the cyclic major index of a cycle permutation $[\pi]$.  Let $[\pi]$ be a cycle permutation of length $n$   and its representative $\hat{\pi}=\hat{\pi}_1\hat{\pi}_2\cdots\hat{\pi}_n$, where $\hat{\pi}_1$ is the largest element in $[\pi]$. The cyclic major index of the cyclic permutation $[\pi]$ is defined to be
\begin{equation}
    {\rm maj}([\pi])={\rm maj}(\hat{\pi}).
\end{equation}
For example, the representative of   the cycle permutation $[4\,1\,3\,2]$ is $\hat{\pi}=4\,1\,3\,2$,  and so its    cyclic major index is defined to be the major index of $\hat{\pi}=4\,1\,3\,2$. It gives that   ${\rm maj}([4\,1\,3\,2])=1+3=4.$

In order to state the cyclic analogue of Stanley's Shuffle Theorem, we will need to  introduce the cyclic descent-bottom  set of a cyclic permutation and recall the splitting map $S_i$ defined by Domagalski, Liang, Minnich, Sagan, Schmidt and Sietsema  in \cite{Domagalski-Liang-Minnich-Sagan-Schmidt-Sietsema-2021},  which maps a cyclic permutation to a linear permutation.    Let $[\pi]$ be a cyclic permutation of length $n$,  the  cyclic descent-bottom set of $[\pi]$ is defined as:
 \begin{equation}\label{defi-cbottoms}
     {\rm cB_d} ([\pi])=\{\pi_{i+1} \colon \pi_i>\pi_{i+1}, \ \text{for} \ 1\leq i\leq n\}
     \end{equation}
   with the convention $\pi_{n+1}=\pi_1$. It should be mentioned that the descent-bottom set of a linear permutation has been
   studied by  Haglund and Visontai \cite{Haglund-Visontai-2012} and Hall and Remmel\cite {Hall-Remmel-2008,Hall-Remmel-2009} respectively.

   It is manifest  from \eqref{defi-cdesnum} and \eqref{defi-cbottoms} that
   \[\#  {\rm cB_d}([\pi])={\rm cdes}([\pi]).\]
   For example,
   \[ {\rm cB_d} ([6\,4\,1\,3])=\{1,4\}.\]

Let $[\pi]$ be a cyclic permutation of length $n$. For $i\in [\pi]$,  Domagalski, Liang, Minnich, Sagan, Schmidt and Sietsema  \cite{Domagalski-Liang-Minnich-Sagan-Schmidt-Sietsema-2021} defined the map $S_i([\pi])$   to be the unique permutation corresponding to  $[\pi]$ which starts with $i$. For example,
\[ S_5([5\,1\,3\,4])=5\,1\,3\,4, \, S_1([5\,1\,3\,4])=1\,3\,4\,5,\, S_3([5\,1\,3\,4])=3\,4\,5\,1,\]
and
\[ S_4([5\,1\,3\,4])=4\,5\,1\,3.\]

We obtain the following generating function of the number of cyclic shufflings of two disjoint cyclic permutations with a given cyclic descent number   and a given cyclic major index.

\begin{thm}[Cyclic Stanley's Shuffle Theorem]\label{cStanley}
Let $[\sigma]\in c\mathcal{S}_m$ and $[\pi] \in c\mathcal{S}_n$ be disjoint cyclic permutations, where ${\rm cdes}([\sigma])=r$ and ${\rm cdes}([\pi])=s$. Moreover, the largest element of $[\sigma]$ and $[\pi]$ is in $[\sigma]$. Then
 \begin{align}\label{cStanley-e}
  & \sum_{[\alpha]\in c\mathcal{S}([\sigma],[\pi]) \atop {\rm cdes}([\alpha])=k}q^{{\rm maj} ([\alpha])}\nonumber\\[5pt]
  &\quad ={m-r+s \brack k-r}{n-s+r-1 \brack k-s-1}q^{{\rm maj} ([\sigma]) +(k-s)(k-r)}
   \sum_{ i \not\in {\rm cB_d}([\pi])}q^{{\rm maj} (S_i([\pi]))}\nonumber\\[5pt]
   &\quad \quad +{m-r+s-1 \brack k-r}{n-s+r \brack k-s}q^{{\rm maj} ([\sigma]) +(k-s+1)(k-r)}
   \sum_{ i \in {\rm cB_d}([\pi])}q^{{\rm maj} (S_i([\pi]))}.
\end{align}
\end{thm}
Setting $q\rightarrow 1$ in Theorem \ref{cStanley},  we  obtain \eqref{AGRR-e}, that is,
\begin{align*}
 & \# c\mathcal{S}([\sigma],[\pi],k) \\[5pt]
  &\quad =\sum_{i \not\in {\rm cB_d}[\pi]}  \dbinom{m-r+s}{k-r}\dbinom{n-s+r-1}{n-k+r}+\sum_{i \in {\rm cB_d}[\pi]} \dbinom{m-r+s-1}{k-r}\dbinom{n-s+r}{n-k+r}\\[8pt]
    &\quad=(n-s)\dbinom{m-r+s}{k-r}\dbinom{n-s+r-1}{n-k+r}+s \dbinom{m-r+s-1}{k-r}\dbinom{n-s+r}{n-k+r}\\[10pt]
    &\quad=\frac{k(m-r)(n-s)+(m+n-k)rs}{(m-r+s)(n-s+r)}\dbinom{m-r+s}{k-r}\dbinom{n-s+r}{k-s}.
\end{align*}

 \section{Proof of Theorem \ref{cStanley}}

This section is denoted to  the proof of  Theorem \ref{cStanley} with the aid of Stanley's Shuffle Theorem \ref{stanley}.

\noindent{\it Proof of Theorem \ref{cStanley}.}  Assume that $[\sigma]\in c\mathcal{S}_m$ and $[\pi] \in c\mathcal{S}_n$ are two disjoint cyclic permutations, where ${\rm cdes}([\sigma])=r$ and ${\rm cdes}([\pi])=s$. Moreover, the largest element of $[\sigma]$ and $[\pi]$ is in $[\sigma]$. Let $\hat{\sigma}=\hat{\sigma}_1\hat{\sigma}_2\cdots \hat{\sigma}_{m}$ be the representative of the cyclic permutation $[\sigma]$, that is,  $\hat{\sigma}_1$ is the largest element of $[\sigma]$. Under the precondition in this theorem, we see that $\hat{\sigma}_1$  is greater than all elements in $[\pi]$. Define
\[\hat{\sigma}'=\hat{\sigma}_2\cdots \hat{\sigma}_{m}.\]
  For $i\in [\pi]$, recall that $S_i([\pi])$ is   the unique permutation corresponding to  $[\pi]$ which starts with $i$.
We claim  that there is a bijection $\psi$ between the set $c\mathcal{S}([\sigma],[\pi])$ and the set $\bigcup_{i\in [\pi]} \mathcal{S}(\hat{\sigma}', S_i([\pi]))$, where  $c\mathcal{S}([\sigma],[\pi])$ denotes the set of cyclic shuffles of $[\sigma]$ and $[\pi]$ and $\mathcal{S}(\hat{\sigma}', S_i([\pi]))$ denotes the set of linear shuffles of $\hat{\sigma}'$ and $S_i([\pi])$.  Moreover,  for $\alpha \in c\mathcal{S}([\sigma],[\pi])$, we have $\psi(\alpha)=\hat{\alpha}'$ such that
\begin{equation}\label{eqn-s1}
 {\rm cdes}([\alpha])= {\rm des}(\hat{\alpha}')+1
\end{equation}
and
\begin{equation}\label{eqn-s2}
 {\rm maj}([\alpha])= {\rm maj}(\hat{\alpha}')+{\rm des}(\hat{\alpha}')+1.
\end{equation}
Let $\alpha \in c\mathcal{S}([\sigma],[\pi])$ and let $\hat{\alpha}=\hat{\alpha}_1\hat{\alpha}_2\cdots \hat{\alpha}_{n+m}$ be the representative of $[\alpha]$, which is a linear permutation corresponding to $[\alpha]$ such that  $\hat{\alpha}_1$ is the largest element in $[\alpha]$. Then $\hat{\alpha}_1=\hat{\sigma}_1$ and
\begin{equation}\label{eqn-t1}
 {\rm cdes}([\alpha])={\rm des}(\hat{\alpha}).
\end{equation}
Define
\[\hat{\alpha}'=\hat{\alpha}_2\hat{\alpha}_3\cdots \hat{\alpha}_{n+m},\]
 which clearly belongs to $\bigcup_{i\in [\pi]} \mathcal{S}(\hat{\sigma}', S_i([\pi]))$. From the construction of $\hat{\alpha}'$ and \eqref{eqn-t1}, we see that  $[\alpha]$ and $\hat{\alpha}'$ satisfy \eqref{eqn-s1} and \eqref{eqn-s2}. Moreover, this process is s reversible. This proved the claim.   Hence
it follows from \eqref{eqn-s1}  and \eqref{eqn-s2}  that
 \begin{align}\label{pf-thm-e3}
  & \sum_{[\alpha]\in c\mathcal{S}([\sigma],[\pi]) \atop {\rm cdes}([\alpha])=k}q^{{\rm maj} ([\alpha])}\nonumber \\[5pt]
   &\quad =\sum_{i\in [\pi]} \sum_{\hat{\alpha}' \in \mathcal{S}(\hat{\sigma}',S_i([\pi]) \atop {\rm des}(\hat{\alpha}')=k-1}q^{{\rm maj} (\hat{\alpha}')+k}\nonumber \\[5pt]
   &\quad =\sum_{ i \not\in {\rm cB_d}([\pi])} \sum_{\hat{\alpha}' \in \mathcal{S}(\hat{\sigma}',S_i([\pi]) \atop {\rm des}(\hat{\alpha}')=k-1}q^{{\rm maj} (\hat{\alpha}')+k}  +\sum_{ i  \in {\rm cB_d}([\pi])}  \sum_{\hat{\alpha}'\in \mathcal{S}(\hat{\sigma}',S_i([\pi]) \atop {\rm des}(\hat{\alpha}')=k-1}q^{{\rm maj} (\hat{\alpha}')+k}.
\end{align}
Observe that ${\rm des}(\hat{\sigma}')={\rm cdes}([\sigma])-1=r-1$ and
 ${\rm des}(S_i([\pi]))={\rm cdes}([\pi])-1=s-1$ if $i \in {\rm cB_d}([\pi])$, otherwise, ${\rm des}(S_i([\pi]))={\rm cdes}([\pi])=s$. Moreover,
 \begin{equation}\label{eqn-s3}
 {\rm maj}(\hat{\sigma}')={\rm maj}([\sigma])-r.
 \end{equation}
  Hence, by Theorem \ref{stanley}, we obtain
 \begin{align}\label{pf-thm-e4}
     &\sum_{ i \not\in {\rm cB_d}([\pi])} \sum_{\hat{\alpha}'\in \mathcal{S}(\hat{\sigma}',S_i([\pi]) \atop {\rm des}(\hat{\alpha}')=k-1}q^{{\rm maj} (\hat{\alpha}')+k} \nonumber \\[5pt]
     &\quad =\sum_{ i \not\in {\rm cB_d}([\pi])} {m-r+s \brack k-r}{n-s+r-1 \brack k-s-1}q^{{\rm maj} (\hat{\sigma}') +{\rm maj} (S_i([\pi]))+(k-s-1)(k-r)+k}
 \end{align}
 and
  \begin{align}\label{pf-thm-e5}
     &\sum_{ i \in {\rm cB_d}([\pi])}  \sum_{\hat{\alpha}'\in \mathcal{S}(\hat{\sigma}',S_i([\pi]) \atop {\rm des}(\hat{\alpha}')=k-1}q^{{\rm maj} (\hat{\alpha}')+k} \nonumber\\[5pt]
     &\quad =\sum_{ i \in {\rm cB_d}([\pi])} {m-r+s-1 \brack k-r}{n-s+r \brack k-s}q^{{\rm maj} (\hat{\sigma}') +{\rm maj} (S_i([\pi]))+(k-s)(k-r)+k} .
 \end{align}
 Substituting \eqref{pf-thm-e4} and \eqref{pf-thm-e5} into \eqref{pf-thm-e3} and using \eqref{eqn-s3}, we obtain \eqref{cStanley-e}. This completes the proof. \qed

 \vskip 0.2cm
\noindent{\bf Acknowledgment.} We are grateful to   Bruce Sagan for bringing this question to our attention and for providing useful comments and suggestions. This work
was supported by   the National Science Foundation of China.

\end{document}